\documentclass[10pt]{article}
 
\usepackage{graphicx}
\usepackage{amssymb}
\usepackage{epstopdf}
\usepackage{amsmath}
\usepackage{amsfonts}
\newcommand{\braket}[2]{\langle #1,#2 \rangle}
\newcommand{\la}{\lambda}

\DeclareSymbolFont{AMSb}{U}{msb}{m}{n}
\DeclareMathSymbol{\N}{\mathbin}{AMSb}{"4E}
\DeclareMathSymbol{\Z}{\mathbin}{AMSb}{"5A}
\DeclareMathSymbol{\R}{\mathbin}{AMSb}{"52}
\DeclareMathSymbol{\Q}{\mathbin}{AMSb}{"51}
\DeclareMathSymbol{\I}{\mathbin}{AMSb}{"49}
\DeclareMathSymbol{\C}{\mathbin}{AMSb}{"43}

\begin{document}
 
\addtolength{\textheight}{0 cm}
\addtolength{\hoffset}{0 cm}
\addtolength{\textwidth}{0 cm}
\addtolength{\voffset}{0 cm}

\setcounter{secnumdepth}{5}
 \newtheorem{proposition}{Proposition}[section]
\newtheorem{theorem}{Theorem}[section]
\newtheorem{lemma}[theorem]{Lemma}
\newtheorem{coro}[theorem]{Corollary}
\newtheorem{remark}[theorem]{Remark}
\newtheorem{ex}[theorem]{Example}
\newtheorem{claim}[theorem]{Claim}
\newtheorem{conj}[theorem]{Conjecture}
\newtheorem{definition}[theorem]{Definition}
\newtheorem{application}{Application}
 
\newtheorem{corollary}[theorem]{Corollary}
\def\HADX{{\cal H}_{\rm AD}(X)}
\def\HADY{{\cal H}_{\rm AD}(Y)}
\def\HADH{{\cal H}_{\rm AD}(H)}
\def\HTADX{{\cal H}_{\rm TAD}(X)}
\def\HTADY{{\cal H}_{\rm TAD}(Y)}
\def\HTADH{{\cal H}_{\rm TAD}(H)}
\def\LX{{\cal L}(X)}
\def\LY{{\cal L}(Y)}
\def\LH{{\cal L}(H)}
 \def\ASD{{\cal L}_{\rm AD}(X)}
 \def\ASDY{{\cal L}_{\rm AD}(Y)}
\def\ASDH{{\cal L}_{\rm AD}(H)}
 \def\ASDP{{\cal L}^{+}_{\rm AD}(X)}
  \def\ASDYP{{\cal L}^{+}_{\rm AD}(Y)}
   \def\ASDHP{{\cal L}^{+}_{\rm AD}(H)}
    \def\TADX{{\cal L}_{\rm TAD}(X)}
        \def\TADY{{\cal L}_{\rm TAD}(Y)}
            \def\TADH{{\cal L}_{\rm TAD}(H)}
 \def\CX{{\cal C}(X)}
\def\CY{{\cal C}(Y)}
\def\CH{{\cal C}(H)}
 
\def\PX{{\cal A}(X)}
\def\PY{{\cal A}(Y)}
\def\PH{{\cal A}(H)}
\def\phi{{\varphi}}
\def\AH{A^{2}_{H}}
\def\B {{\cal B}} 
\def\C {{\cal C}} 
\def\H{{\cal H}}

\newcommand{\al}{\alpha}
\newcommand{\de}{\delta}

\newcommand{\ra}{\rightarrow}
\def\phi{{\varphi}}

\title{A variational theory for monotone vector fields}
\author{ Nassif  Ghoussoub\thanks{Research partially supported by a grant from the Natural Sciences and Engineering Research Council of Canada. } 
\\
\small Department of Mathematics,
\small University of British Columbia, \\
\small Vancouver BC Canada V6T 1Z2 \\
\small {\tt nassif@math.ubc.ca}\\
{\it   Dedicated to Felix Browder on his 80th birthday}
}
\maketitle

\begin{abstract} Monotone vector fields were introduced almost 40 years  ago as  nonlinear extensions of positive definite linear operators, but also as natural  extensions of gradients of convex potentials.   These vector fields are not always derived from potentials in the classical sense, and as such they are not always amenable to the standard methods of the calculus of variations.   We describe here how the selfdual variational calculus developed recently by the author, provides a  variational approach to PDEs and evolution equations driven by maximal monotone operators.  To any such a vector field $T$ on a reflexive Banach space $X$, one can associate a convex selfdual Lagrangian $L_T$ on phase space $X\times X^*$ that can be seen as a  ``potential" for $T$, in the sense that the problem of inverting  $T$ reduces to  minimizing the convex energy $L_T$. This variational approach to maximal monotone operators allows their theory to be analyzed with the full range of methods --computational or not-- that are available for variational settings.  Standard convex analysis (on phase space) can then be used to establish many old and new results concerned with the identification, superposition, and resolution of such vector fields.

\end{abstract}

\section{Introduction} Monotone vector fields are 
those --possibly  set valued and nonlinear-- operators $T$ from a Banach space $X$ into (the subsets of) its dual $X^*$, whose graphs $G(T)=\{ (x,p)\in X\times X^*; p\in T(x)\}$ are {\it monotone} subsets of 
$X\times X^*$, i.e., they satisfy:
\begin{equation}
\hbox{$\langle u-v, p - q \rangle \geq 0$ for every $(u,p)$ and $(v, q)$ in $G(T)$.}
\end{equation}
The {\it effective domain} $D(T)$ of $T$ is then the set of all $u\in X$ such that $T(u)$ is nonempty. 
A useful  subclass consists of the so-called  {\it maximal monotone operators} which refer to those monotone operators whose graph $G(T)$ is maximal in the family of monotone subsets of $X\times X^*$, ordered by set inclusion. Starting in the sixties, this theory was studied in depth by G.J. Minty, F. Browder, H. Brezis, L. Nirenberg, and T. Rockafellar to name a few.  Felix Browder was one of the pioneers in developing a systematic approach to study maximal monotone operators and their role in connection with nonlinear partial differential equations and other aspects of nonlinear analysis. See for example  \cite{Bro1}, \cite{Bro2} and his other  numerous  contributions to this subject as referenced in the monographs of Brezis \cite{Br},    Phelps \cite{Ph}, and  Kindehlehrer-Stampachia \cite{KS}.

Now it is well known that many basic linear and nonlinear elliptic PDEs are variational and can be written in the form 
\begin{equation}
\hbox{$\partial \Phi (u)=p$ on $H$ \quad or \quad  $-{\rm div}(\partial \phi (\nabla u(x))+\lambda u(x)=p(x)$ on $\Omega\subset \R^n$,}
\end{equation}
 where $\Phi$ (resp., $\phi$) is a convex functional on an infinite dimensional function space $H$ (resp., $\R^n$). They are the Euler-Lagrangian equations associated to the energy functional 
\begin{equation}
 \hbox{$I(u)=\Phi (u)-\langle u, p\rangle$ \quad resp., \quad $J(u)=\int_\Omega \big\{\phi (\nabla u (x))+\frac{\lambda}{2} |u(x)|^2- u(x) p(x) \big\}dx$.}
\end{equation}
 However,  a large number of PDEs can be formulated as
\begin{equation}\label{one}
\hbox{$Tu=p$ \quad or \quad $-{\rm div}(T(\nabla u(x))+\lambda u(x)=p(x)$ on $\Omega\subset \R^n$,}
\end{equation}
where $T$ is a vector field on $H$ (resp., $\R^n$) that is not derived from a potential, yet it shares many properties with $\partial \Phi$ (resp., $\partial \phi$) such as  monotonicity.
Their solutions cannot therefore be obtained by the classical methods of the calculus of variations, leading people to resort to nonvariational techniques as described by  C. Evans in  chapter 9 of   his  landmark book \cite{Ev}.

Moreover, initial-value problems of the form $u(0)=u_0$
\begin{equation}\label{two}
 \hbox{$\dot u (t)+T (u(t))=p(t)$ \quad or \quad $\dot u (t)-{\rm div}(T (\nabla u(t))=p(t)$\quad on $[0,1]$} 
\end{equation}
are not variational in the classical sense, even when $T$ is a potential operator (i.e., $T=\partial \phi$). 

Our goal in this  paper is to show that there is indeed a variational theory for these equations. More precisely,  we show  that one can associate to the vector field $T$ a convex Lagrangian $L_T$ on phase space $H\times H$ (resp., on $\R^n\times \R^n$) such that the equations  in $(\ref{one})$ can be solved  by simply minimizing the functional 
\begin{equation}
\hbox{$I(u)=L_T(u, p) -\langle u, p\rangle$\quad  resp., \quad  $J(u)=\int_\Omega \big\{L_T\big(\nabla (-\Delta)^{-1}(-\lambda u +p), \nabla u\big) +\lambda |u|^2-  up \big\}dx$,}
\end{equation}
on $H$ (resp., $H^1_0(\Omega)$), where here $-\Delta$ is considered as an operator from $H^1_0(\Omega)$ onto $H^{-1}(\Omega)$. 
Similarly, we shall be able to solve the above evolution equations  on a time interval $[0, 1]$ say, by minimizing
\begin{equation}
{\cal I}(u)=\int_0^1\left\{L_T(u(t), -\dot u (t)+p(t)) -\langle u(t), p(t) \rangle \right\}dt +\ell_{u_0} (u(0), u(1)),
\end{equation}
respectively
\begin{equation}
{\cal J}(u)=\int_0^1\int_\Omega \left\{L_T\big( \nabla (-\Delta)^{-1}(-\dot u +p), \nabla u\big) -  u p  \right\}dx dt +\ell_{u_0} (u(0), u(1)),
\end{equation}
on $W^{1,2}([0, 1]; H)$ (resp., on $W^{1,2}([0, 1]; H^1_0(\Omega))$), and  where $\ell_{u_0}$ is another convex Lagrangian associated to the given initial condition $u(0)=u_0$.  
 
Our approach was motivated indirectly by a 1975 conjecture of Brezis-Ekeland \cite{BE} which eventually led us to develop in a series of papers (\cite{G2}--\cite{GT2}) the concept and the calculus of selfdual Lagrangians $L$ on phase space $X\times X^*$ where $X$ is a reflexive Banach space. This allowed us  to provide variational formulations and resolutions for various differential equations which are not variational in the sense of Euler-Lagrange. The main idea behind this theory --which is summarized in the upcoming monograph \cite{G10}--  originated from the fact that a large set of PDEs and evolution equations can be written in the form 
\begin{equation} \label{basic.equation}
(p,u)\in \partial L(u,p),
\end{equation}
where $\partial L$ is the subdifferential of a   Lagrangian $L:X\times X^*\to \R\cup \{+\infty\}$ that is convex  and lower semi-continuous --in both variables-- while satisfying the {\it selfdual} conditions,
\begin{equation}
\hbox{$L^*( p, u) =L(u, p )$   for all $(p,u)\in X^{*}\times X$.}
\end{equation}
Here $L^*$ is the Legendre transform of $L$ in both variables, that is
\[
L^*(p,u)=\sup\left\{ \langle p,y\rangle +\langle u,q\rangle -L(y, q);\, (y, q)\in X\times X^*\right\}.             
\]
As will be shown in Lemma 2.1 below, solutions can then be found for a given $p$,  by simply minimizing the functional $I_p(u)=L(u,p)-\langle u,p\rangle$ and by proving that the infimum is actually zero. In other words, by defining the  following  vector fields of L at $u\in X$ 
 to be the --possibly empty-- sets
\begin{equation} \label{selfdual.vector}
\hbox{${\overline\partial} L (u):=\{p\in X^*; \, L(u,p)-\langle u,p\rangle =0\}=\{p\in X^*; \, (p,u)\in \partial L(u,p)\}$,}
\end{equation}
one can then find variationally the zeroes of  those set-valued maps $T: X\to 2^{X^*}$ of the form $T(u)=\bar\partial L(u)$ for some selfdual Lagrangian $L$ on $X\times X^*$, even when such maps are not derived from ``true potentials". 

These {\it selfdual vector fields} are also 
natural extensions of  subdifferentials of convex lower semi-continuous functions, in which case the corresponding  selfdual Lagrangians is $L(u,p)=\phi (u) +\phi^*(p)$ where $\phi$ is  such a  function on $X$, and $\varphi^{*}$ is its Legendre conjugate on $X^{*}$. It is then easy to see that ${\overline \partial}  L (u)=\partial \phi (u)$, and that 
 the corresponding variational problem (i.e., minimizing $I(u)=L(u,p)=\phi (u) +\phi^*(p)-\langle u, p\rangle $) reduces to the classical approach of minimizing a convex functional in order to solve equations of the form $p\in \partial \phi (u)$. 

More interesting examples of selfdual Lagrangians are of the form 
$
L(u,p)=\varphi (u) +\varphi^{*}(-\Gamma u+p)
$
 where $\varphi$ is a convex and lower semi-continuous function on $X$, and  $\Gamma:X\to X^{*}$ is a skew-symmetric operator. The corresponding selfdual vector field is then, 
\[
\hbox{${\overline \partial} L (u)=\Gamma u+\partial \phi (u)$.}
\]
More generally,  if the operator $\Gamma$ is merely non-negative (i.e., $\langle \Gamma u, u\rangle \geq 0$), then one can still write the vector field  $\Gamma +\partial \phi$ as ${\overline \partial} L$ for some selfdual Lagrangian $L$ defined now on $X\times X^*$, as
\[
L(u,p)=\psi (u) +\psi^{*}(-\Gamma^{as} u+p)
\]
where $\psi$ is the convex function $\psi (u)=\frac{1}{2} \langle \Gamma u, u\rangle +\varphi (u)$ and $\Gamma^{as}=\frac{1}{2}(\Gamma -\Gamma^{*})$ is the anti-symmetric part of $\Gamma$, and $\Gamma^{sym}=\frac{1}{2}(\Gamma +\Gamma^{*})$ is its symmetric part. The main interest being that equations of the form  $p\in \Gamma u+\partial \phi(u) $ can now be solved for a given  $p\in X^*$, by simply minimizing the functional 
\[
I_p(u)=\psi (u) +\psi^{*}(\Gamma^{as} u+p)-\langle u, p\rangle
\]
 and proving that its infimum is actually zero.

It was therefore natural to investigate the relationship between maximal monotone operators and seldual vector fields since both could be seen as extensions of the superposition of subgradients of convex functions with skew-symmetric operators. An early indication was the observation we made in \cite{G2}, that selfdual vector fields are necessarily maximal monotone. We suggested calling them  then  {\it ``integrable maximal monotone fields"}  not suspecting that one could eventually prove that all maximal monotone operators are integrable in the sense that they all do derive from selfdual Lagrangians  \cite{G5}. This surprising development actually occured when we realized through the book of Phelps \cite{Ph} that  
Krauss \cite{K} and Fitzpatrick \cite{F}  had done some work in this direction in the 80's, and had managed to associate to a maximal monotone operator $T$, a ``sub-selfdual Lagrangian", i.e., a convex lower semi-continuous function $L$ on state space $X\times X^*$ satisfying 
\begin{equation}\label{sub}
\hbox{$L^*(p,u) \geq L(u,p) \geq \langle p, u\rangle$ on $ X\times X^*$,}
\end{equation}
 in such a way that (in our terminology) $T=\bar \partial L$.   The question whether one can  establish the existence of a truly selfdual  Lagrangian associated to $T$, was actually one of the original questions of Kirkpatrick  \cite{F}. We eventually stated the following result in \cite{G5} whose complete proof is included in section 1.

 \begin{theorem}  \label{ghoussoub} Let $L$ be a proper selfdual Lagrangian $L$ on a reflexive Banach space $X\times X^*$, then the vector field $u\to \bar \partial L (u)$ is maximal monotone.

 Conversely, if $T: D(T) \subset X\to 2^{X^*}$ is a maximal monotone operator with a non-empty domain, then there exists a selfdual Lagrangian $L$ on $X\times X^*$ such that $T= \bar \partial L$. 
 
  \end{theorem}
  
 Almost one year later, we eventually learned that the sufficient condition in Theorem \ref{ghoussoub} had been  established by R.S. Burachik and B. F. Svaiter in \cite{BS}, while the necessary condition was shown by B. F. Svaiter in \cite{S}. The methods in both directions are quite different  from those described here. Actually, in  our original method for the necessary condition, we used Asplund's averaging technique between the sub-selfdual Lagrangian $L$ given by Fitzpatrick, and its Legendre dual $L^*$. This turned out to warrant additional boundedness assumption which required an additional approximation argument. However,  upon seeing our paper, Baushke and Wang \cite{BW}  noted that ``proximinal interpolaton" between  $L$ and $L^*$  gives an explicit formula for the selfdual Lagrangian. It is this formula that we adopt here since --unlike Svaiter's proof that  relies on Zorn's lemma-- this formula insures a measurable selection of selfdual Lagrangians to correspond to a  measurable family of  maximal monotone operators. This is  an important technical issue in the study of evolution equations driven by  time-dependent vector fields. 

It is  worth comparing at this stage, the above result with the following celebrated result of Rockafellar \cite{Ph}, which gives an integral representation of those maximal monotone operators $T: X\to 2^{X^*}$ that are {\it cyclically monotone}, i.e., those that verify for any finite number of points $(u_i, p_i)_{i=0}^n$ in the graph $G(T)$ with $u_0=u_n$, we have 
\begin{equation}
\hbox{$\sum\limits_{i=1}^n\langle p_k, u_k-u_{k-1} \rangle \geq 0$.}
\end{equation}
\begin{theorem} [Rockafellar] If $\phi:X\to \mathbb{R}\cup\{+\infty\}$ is a proper convex and lower semi-continuous functional on a Banach space $X$, then its subdifferential map  $u\to \partial \phi (u)$ is a maximal cyclically monotone map.

 Conversely if $T:X\to 2^{X^*}$ is a maximal cyclically monotone map with a non-empty domain, then there exists a a proper convex and lower semi-continuous functional on  $X$ such that $T=\partial \phi$.
 \end{theorem}
In this paper, we shall emphasize the advantages to identifying  maximal monotone operators with selfdual vector fields. Here are some of them.\begin{trivlist}
\item $\bullet$\, As mentioned above, all equations, systems, variational inequalities, and  dissipative initial value  parabolic problems which traditionally involve maximal monotone operators,  can now be formulated and resolved variationally.   These problems can therefore be analyzed with the full range of methods --computational or not-- that are available for variational settings.  We shall describe some concrete examples in section 6. 

\item $\bullet$\,   While issues around the various ways to combine maximal monotone operators are often delicate to prove,  the class of selfdual Lagrangians  possesses remarkable  permanence properties that are also easy to establish.   Indeed, operations such as superposition, direct sum and convolution of Lagrangians reduce to standard convex analysis on phase space,  making the calculus of selfdual Lagrangians (and therefore of maximal monotone operators) as manageable as convex analysis, yet much more encompassing  \cite{G2}. In section 3, we shall describe  the Lagrangian calculus that correspond to the various operations on maximal monotone operators. 

\item $\bullet$\, Selfduality allows for the superposition of appropriate boundary Lagrangians  with ``interior"   Lagrangians,  leading to the resolution of problems with various linear and nonlinear boundary constraints that are not amenable to standard variational theory \cite{G4}, \cite{GT2}. We shall describe this aspect in section 5, by concentrating on time boundary conditions in evolution equations. 

\item $\bullet$\, Selfdual Lagrangians defined on state spaces ``lift" to selfdual Lagrangians on path spaces leading to a unified approach for stationary and dynamic equations.  More precisely,  flows of the form  $-\dot u (t) \in T(t, u(t))$ corresponding to time-dependent maximal monotone operators with a variety of time-boundary conditions can be reformulated and resolved  as $0\in \overline \partial {\cal L} (u)$ where ${\cal L}$ is a corresponding selfdual Lagrangian on path space, a phenomenon that leads to  natural and quite interesting iterations (See \cite{G2}, \cite{GT3}).  We shall also recover variationally --in section 5-- the construction of a semi-group of contractions associated to a maximal monotone operator. 

\item $\bullet$\,  The class of anti-symmetric Hamiltonians that one can associate to selfdual Lagrangians (\cite{G3}, \cite{GM4},  \cite{GM2},  \cite{GM3}, \cite{GM4}) goes beyond the theory of maximal monotone operators, and  leads to a much wider array of applications. It shows among other things that they can be superposed with certain nonlinear operators that are far from being maximal monotone \cite{G3}, \cite{G4} and \cite{G6}. The corresponding  class of PDEs --including Hamiltonian systems, Cauchy-Riemann and Navier-Stokes equations-- and their variational principles will not be considered here, but are studied in detail in the upcoming monograph \cite{G10}.  

  \end{trivlist}
 
 \section{Selfdual vector fields are maximal monotone operators}
 
 To any convex lower semi-continuous Lagrangian $L:X\times X^* \to \R \cup\{+\infty\}$, we can associate  a vector field $ \delta L$ at any $x\in X$,  
 to be the --possibly empty-- set
\begin{equation} \label{vector,field}
\hbox{${\delta} L (x):=\{p\in X^*; \, (p,x)\in \partial L(x,p)\}$.}
\end{equation}
It is easy to see that the convexity of $L$ yields that $x\to \delta L (x)$ is a monotone map, since if $p\in \delta L(x)$ and $q\in \delta L(y)$, then it can be easily verified that
\[
\langle p-q, x-y\rangle =\frac{1}{2}\langle (p,x)-(q,y), (x,p)-(y,q) \rangle \geq 0.
\]
We can also associate to $L$ another --not  necessarily monotone-- vector field as: 
\begin{equation} \label{selfdual.vector}
\hbox{${\overline\partial} L (x):=\{p\in X^*; \, L(x,p)-\langle x,p\rangle =0\}.$}
\end{equation}

As the following lemma indicates,  $\bar \partial L$ should not be confused with the subdifferential $\partial L$ of $L$ as a convex function on $X\times X^*$. 

\begin{lemma} Let  $L$ be a proper convex lower semi-continuous Lagrangian on $X\times X^*$, where $X$ is a reflexive Banach space.
\begin{enumerate}
\item If $L$   satisfies 
$\hbox{$L(x,p)\geq \langle x, p\rangle$ for all $(x,p)\in X\times X^*$}$,
then $\bar \partial L (x) \subset \delta L (x)$.

\item If $L$ also satisfies
$\hbox{$L^*(p,x)\geq  L(x,p)$ for all $(x,p)\in X\times X^*$}$,
then $\bar \partial L (x)=\delta L (x)$.
\item  In particular, if $L$ is a selfdual Lagrangian on $X\times X^*$,  then  $\bar \partial L (x)=\delta L (x)$.
\end{enumerate}
\end{lemma}
\noindent{\bf Proof:} 1) Assuming $p\in {\overline \partial}  L (x)$, we consider any $(y, q)\in X\times X^*$ and write
\begin{eqnarray*}
L(x+y, p+q) -L(x,p)&\geq &t^{-1} \big[L(x+ty, p+tq)-L(x,p)\big]\\
&\geq&t^{-1} \big[\langle x+ty, p+tq\rangle -\langle x,p\rangle \big]\\
&\geq&\langle x, q\rangle +\langle y, p\rangle +t\langle y,q\rangle
\end{eqnarray*}
Letting $t\to 0^+$, we get that $L(x+y, p+q) -L(x,p)\geq \langle x, q\rangle +\langle y,p\rangle$ which means that we have $(p, x)\in \partial L(x,p)$, and consequently $p\in \delta L (x)$. 

 2)\,  Indeed, assume $p\in \delta L (x)$ which means that $(p,x)\in \partial L(x,p)$ and therefore by  Legendre-Fenchel duality applied to  $L$ on $X\times X^*$, we have
\[
 L^*(p, x)+L(x, p) = 2\langle x, p\rangle.  
 \]
 But since $L^*(p, x)\geq L(x,p)\geq \langle x, p\rangle$, we must have that  $L^*(p, x)-\langle x, p\rangle =  L(x,p) - \langle x, p\rangle =0$, and therefore $p\in \bar \partial L (x)$. 
 
 3)\, It suffices to show that a selfdual Lagrangian $L$ satisfies $L(x,p) \geq \langle x, p\rangle$ for all $(x,p)\in X\times X^*$. But this follows easily from Legendre-Fenchel duality applied to $L$ since then
 \[
 2L^*(p,x)=2L(x,p)= L^*(p, x)+L(x, p) \geq  2\langle x, p\rangle.  
  \]
    \hfill $\square$.

For the sequel, we note that  $p\in {\overline \partial}  L (x)$ if and only if  $0\in {\overline \partial}  L_p (x)$ where $L_p$ is the selfdual Lagrangian $L_p(x, q)=L(x, p+q)+\langle x, p\rangle$. This is also equivalent  to the statement that the infimum of the functional $I_p(u)=L(u, p)-\langle u, p\rangle$ is zero and is attained at $x\in X$. This leads to the following proposition which is a particular case of a more general result established in \cite{G2}.

\begin{proposition} \label{main.existence} Let $L$ be a convex lower semi-continuous  selfdual Lagrangian  on a reflexive Banach space $X\times X^{*}$, such that  for some $x_0\in X$, the function $p\to L(x_0,p)$ is  bounded on the balls of $X^*$.   Then for each $p\in X^*$, there exists $\bar x\in X$ such that: 
 \begin{equation}
 \left\{ \begin{array}{lcl}
\label{eqn:main.existence}
  L( \bar x, p)-\langle \bar x, p\rangle&=&\inf\limits_{x\in X}\left\{ L(x,p)-\langle x, p\rangle\right\}=0.\\
 \hfill  p &\in & \bar \partial L (\bar x).
\end{array}\right.
 \end{equation}
 \end{proposition}
\noindent{\bf Proof:} We can assume that $p=0$ by considering the translated  Lagrangian  $M(x, q)=L(x, p+q)-\langle x, p\rangle$ which is also selfdual on $X\times X^*$. In this case $M(x,0) \geq 0$ for every $x\in X$.

 Now consider  $({\cal P}_{q})$ the primal  minimization problem
$h(q)=\inf\limits_{x\in X}M(x,q)$
in such a way that $({\cal P}_{0})$ is the initial problem $h(0)=\inf\limits_{x\in X}M(x,0)$, 
and the dual problem $({\cal P}^{*})$ is therefore
$\sup\limits_{y\in X}-M^{*}(0,y)$ (See for example \cite{ET}). 

We readily have the following weak duality formula:
\[
\inf {\cal P}_{0}:=\inf_{x\in X}M(x,0) \geq 0\geq \sup_{y\in X}-M^{*}(0,y):=\sup{\cal P}^{*}. 
\]
Note that  $h$ is convex on $X^{*}$, and that its Legendre conjugate satisfies  for all $y\in X$.
\[
h^{*}(y)=M^{*}(0,y)=M(y,0)=L(y, p)-\langle y, p\rangle
\]
Moreover,  $h(q)=\inf\limits_{x\in X}M(x,q) \leq L(x_0, p+q)-\langle x_0, p\rangle$ and therefore $q\to h(q)$   is bounded above on the balls of  $X^{*}$, 
and hence    
it is subdifferentiable at $0$ (i.e., the problem $({\cal P}_{0})$ is then stable). Any point $\bar x \in \partial h (0)$ satisfies $h(0) +h^{*}(\bar x)=0$, which means that
\[
-\inf_{x\in X}M(x,0)=-h(0)=h^{*}(\bar x)=M^{*}(0, \bar x)=M(\bar x,0)\geq  \inf_{x\in X}M(x,0).
\]
It follows that $\inf_{x\in X}M(x,0)=M(\bar x,0)\leq 0$ and the infimum of $({\cal P})$  is therefore zero and is attained at $\bar x$.  

\begin{lemma} \label{inf.convolution} Let $L$ be a proper selfdual Lagrangian on $X\times X^*$, then for any  convex continuous function $\phi$ on $X$, the Lagrangian defined by
\[
M(x,p)=\inf\left\{L(x, p-r)+\phi (x) +\phi^*(r); \, r\in X^*\right\}
\]
is also selfdual on $X\times X^*$.
\end{lemma}
 {\bf Proof:} Indeed,  fix $(q, y)\in X^{*}\times X$ and write:
\begin{eqnarray*}
 M^{*} (q,y)
  &=&\sup\{\langle q, x\rangle +  \langle y, p\rangle- L(x, p-r)- \phi (x) - \phi^*(r); (x, p,r)\in X\times X^*\times X^{*}\}\\
&=&  \sup\{\langle q, x\rangle +  \langle y, r+s\rangle- L(x, s)- \phi (x) - \phi^*(r); (x, s, r)\in X\times X^*\times X^{*}\}\\
&=&\sup_{x\in X} \left\{\langle x,q\rangle +\sup\limits_{(s, r)\in X^*\times X^*}\{\langle y, r+s\rangle -L(x,s)- \phi^*(r)\}- \phi (x)\right\}\\
&=&\sup_{x\in X}\left\{\langle x,q\rangle +\sup_{s\in X^*}\{\langle y,s\rangle -L(x,s)\} +\sup_{r\in X^*} \{\langle y,r\rangle -\phi^* (r)\}-\phi (x)\right\}\\
&=&\sup_{x\in X}\left\{\langle x,q\rangle +\sup_{s\in X^*}\{\langle y,s\rangle -L(x,s)\} + \phi (y)-\phi (x)\right\}\\
&=&\sup_{x\in X}\sup_{s\in X^*}\left\{\langle x,q\rangle +\langle y,s\rangle  -L(x,s)-\phi (x)\right\} + \phi (y)\\
&=& (L+T_\phi)^{*}(q,y)+\phi (y)
\end{eqnarray*}
where $T_\phi (x,s):=\phi (x) $ for all $(x,s)\in X\times X^{*}$. Note now that 
\begin{eqnarray*}
T_\phi^*(q,y)=\sup_{x,s}\left\{\langle q,x\rangle +\langle y,s\rangle -  \phi (x)\right\}=\left\{\begin{array}{lll}+\infty &\hbox{if }&y\ne 0\\ \phi^*(q)&\hbox{if }&y=0\end{array} \right.
\end{eqnarray*}
in such a way that by using the duality between sums and inf-convolutions in both variables, we get
\begin{eqnarray*}
(L+T_\phi)^{*}(q,y)&=& L^*\star T_\phi^*\,(q,y)\\
&=&\inf_{r\in X^{*},z\in X}\left\{ L^*(r,z)+T_\phi^*(-r+q,-z+y)\right\}\\
&=&\inf_{r\in X^*}\left\{ L^*(r,y)+ \phi^*(-r+q)\right\}
\end{eqnarray*}
and  finally 
\begin{eqnarray*}
M^{*} (q,y)&=&(L+T_\phi)^{*}(q,y)+\phi (y)\\
&=&\inf_{r\in X^*}\left\{ L^*(r,y)+\phi^*(-r+q)\right\}+\phi (y))\\
&=&\inf_{s\in X^*}\left\{ L^*(q-s,y)+\phi^*(s)\right\}+\phi (y))\\
&=&\inf_{s\in X^*}\left\{ L(y,q-s)+\phi (y) +\phi^*(s) \right\}\\
  &=& M(q,y).
\end{eqnarray*}

 \begin{proposition} \label{subdifferential} Let $L$ be a selfdual Lagrangian $L$ on a reflexive Banach space $X\times X^*$. The following assertions then hold:
  \begin{enumerate}
  \item The vector field $x\to \bar \partial L (x)$ is maximal monotone.
 
  \item  If $L$ is strictly convex in the second variable, then the maximal monotone vector field  $x\to \bar \partial L (x)$ is single-valued on its domain.
   
   \item If  $L$ is uniformly convex in the second variable (i.e.,  if  
   $L(x,p)-\epsilon\frac{\| p\|^2}{2}$  is convex in $p$ for some $\epsilon>0$), 
    then the vector field  $x\to \bar \partial L (x)$ is a Lipschitz  maximal monotone operator on its domain. 
  \end{enumerate}
    \end{proposition}   
{\bf Proof:} Denoting by $J:X\to 2^{X^*}$  the duality map between $X$ and $X^*$, that is
\[
J(x)=\{p\in X^*;\, \langle x, p\rangle =\|x\|^2\}.
\]
In order to show that $\bar\partial L$ is maximal monotone,  it suffices to show that the vector field $\bar\partial  L+    J$ is onto \cite{Ph}.   In other words, we need to find for any given $p\in X^*$,  an $x\in  X$ such that $p\in \bar\partial  L(x)+  J(x)$.
For that, we consider the following Lagrangian on $X\times X^*$.
\[
M(x,p)=\inf\left\{L(x, p-r)+\frac{1}{2}\|x\|_X^2 + \frac{1}{2}\|r\|_{X^*}^2; \, r\in X^*\right\}.
\]
It is a selfdual Lagrangian according to the previous lemma. Moreover -- assuming without loss of generality-- that the point $(0,0)$ is in the domain of $L$, we get   the estimate 
$
  M(0,p)\leq L(0,0)+\frac{1}{2}\|p\|_{X^*}^2  
$, and therefore Proposition \ref{main.existence} applies and we obtain $\bar x \in X$ so that $ p\in \bar\partial  M ( \bar x)$.
This means that 
\[
M(\bar x,p)-\langle \bar x, p\rangle=\inf\left\{L(\bar x, p-r)-\langle \bar x, p-r\rangle+\frac{1}{2}\|\bar x\|_X^2 + \frac{1}{2}\|r\|_{X^*}^2-\langle \bar x, r\rangle;\,  r\in X^* \right\}=0, 
\]
which means that there exists $\bar r \in X^*$ such that
\[
\hbox{$L(\bar x, p-\bar r)-\langle \bar x, p-\bar r\rangle=0$ \quad and \quad $\frac{1}{2}\|\bar x\|_X^2 + \frac{1}{2}\|\bar r\|_{X^*}^2-\langle \bar x, \bar r\rangle=0$.}
\]
In other words, there exists $\bar r\in  J(\bar x)$ such that  $p - \bar r  \in  \bar\partial  L ( \bar x)$ and we are done.  

The other assertions of the proposition are straightforward and left to the interested reader. 
 
\section{Maximal monotone operators are selfdual vector fields}

We start with the following lemma which is essentially due to Fitzpatrick \cite{F}. 

\begin{lemma} \label{fitz}Let $T: D(T) \subset X\to 2^{X^*}$ be a monotone operator, and consider  on $X\times X^*$ the Lagrangian $L_T$ defined by
\begin{equation}
L_T(x,p)=\sup\{\langle p,y\rangle +\langle q, x-y\rangle; \, (y,q)\in G(T)\}
\end{equation}
\begin{enumerate}
\item If $D(T)\neq \emptyset$, then $L_T$ is a convex and lower semi-continuous function on $X\times X^*$ such that for every $x\in D(T)$, we have $Tx\subset \bar \partial L_T(x) \cap \delta L_T (x)$. Moreover, we have 
\begin{equation}
\hbox{$L_T^*(p,x) \geq L_T(x,p)$ for every $(x,p)\in X\times X^*$.}
\end{equation}

\item If $T$ is maximal monotone, then $T= \bar \partial L_T=\delta L_T$
and 
\begin{equation} \label{fenchel}
\hbox{$L_T(x,p)\geq \langle x,p\rangle$ for all $(x,p)\in X\times X^*$.} 
\end{equation}
\end{enumerate}
\end{lemma} 

{\bf Proof:} (1) If $x\in D (T)$ and  $p\in Tx$, then the monotonicity of $T$ yields for any $(y, q)\in G(T)$
\[
\langle x,p\rangle \geq \langle y,p\rangle +\langle x-y, q\rangle
\]
in such a way that $L_T(x,p)\leq \langle x,p\rangle$. On the other hand, we have 
\[
L_T(x,p) \geq \langle x,p\rangle +\langle p, x-x\rangle =\langle x,p\rangle,
\]
and therefore $p\in \bar \partial L_T(x)$.

Write now for any $(y, q)\in X\times X^*$, 
\begin{eqnarray*}
L_T(x+y, p+q)-L_T(x,p)&=&\sup\left\{ \langle p+q, z\rangle +\langle r, x+y\rangle-\langle z,r\rangle;\, (z,r)\in G(T)\right\}-L_T(x,p)\\
&\geq &\langle p+q, x\rangle +\langle p, x+y\rangle-\langle p,x\rangle -\langle p,x\rangle\\
&= &\langle q, x\rangle +\langle p, y\rangle 
\end{eqnarray*}
which means that $(p,x) \in \partial L_T (x,p)$ and therefore $p\in \delta L_T(x)$.

Note also that $L_T(x,p)=M_T^*(p,x)$ where $M_T$ is the Lagrangian on $X\times X^*$ defined by 
\begin{equation} \label{hamiltonian}
 M_T(x,p)= \left\{ \begin{array}{lcl}
\hbox{$\langle x, p \rangle$ \quad if $(x,p)\in G(T)$}\\
\hbox{$+\infty$ \quad \quad \quad otherwise.} 
\end{array}\right.
\end{equation}
Since $L_T(x, p)=\langle x, p \rangle=M_T(x,p)$ whenever $(x,p)\in G(T)$, it follows that 
$L_T\leq M_T$ on $X\times X^*$ and so $L_T^*(p,x)\geq M_T^*(p,x)=L_T(x,p)$ everywhere. \\

(2) If now $T$ is maximal then necessarily $Tx= \delta L_T(x)\cap \bar \partial L_T(x)=\delta L_T(x)$ since $x\to \delta L_T(x)$ is a monotone extension of $T$.  

In order to show (\ref{fenchel}), assume to the contrary that $L_T(x,p)< \langle x,p\rangle$ for some $(x,p)\in X\times X^*$. It follows that 
\[
\hbox{$\langle p, y\rangle +\langle q, x-y\rangle < \langle p, x\rangle$ for all $(y,q)\in G(T),$}
\]
and therefore 
\[
\hbox{$\langle p-q, x-y\rangle> 0$ for all $(y,q)\in G(T).$}
\] 
But since $T$ is maximal monotone, this means that $p\in Tx$. But then $p\in \bar \partial L_T(x)$ by the first part, leading to $L_T(x,p)= \langle x,p\rangle$, which is a contradiction.

Finally, note that property (\ref{fenchel}) on $L_T$ yields that  $\bar \partial L_T(x) \subset \delta L_T (x)$ and therefore  $Tx=\bar \partial L_T(x)$.

\begin{proposition} Let  $X$ be a separable reflexive Banach space, and let $L$ be a  convex lower semi-continuous Lagrangian on $X\times X^*$ that satisfies 
 \begin{equation}
\hbox{$L^*(p,x) \geq L(x,p)\geq \langle x, p\rangle$ for every $(x,p)\in X\times X^*$.}
\end{equation}
Then, there exists a selfdual Lagrangian $N$ on $X\times X^*$ such that $\bar \partial L=\bar \partial N$ and 
 \begin{equation}
\hbox{$ L(x,p) \leq N(x,p) \leq L^*(p,x)$ for every $(x,p)\in X\times X^*$.}
 \end{equation}
\end{proposition} 

{\bf Proof:}   The Lagrangian $N$ is simply  the {\it proximal average}  between  $L$ and $\tilde L$ where $\tilde L (x, p)=L^*(p, x)$. It is defined as
 \[
 N(x, p):=\inf \left\{\frac{1}{2}L(x_1, p_1)+\frac{1}{2}L^*(p_2, x_2)+\frac{1}{8}\|x_1-x_2\|^2+\frac{1}{8}\|p_1-p_2\|^2; \, (x, p)=\frac{1}{2}(x_1, p_1) + \frac{1}{2}(x_2, p_2)\right\}.
 \]
It is easy to see that  $L(x, p) \leq N(x,p) \leq L^*(p,x)$.  Before showing that it is a selfdual Lagrangian on $X\times X^*$,  we note that $ \bar \partial L(x)= \bar \partial N (x)$. Indeed, first  it is  clear that  $ \bar \partial N (x)\subset \bar \partial L(x)$. On the other hand, since $L(x,p)\geq \langle x, p\rangle$, we have from Lemma 2.1  that $ \bar \partial L(x)\subset \delta  L(x)$  which means that if $p\in  \bar \partial L(x)$ then  $(p,x)\in \partial L(x,p)$ and therefore
$
L(x,p)+L^*(p,x)=2\langle x,p\rangle.
$
Again, since  $ p\in \bar \partial L(x)$  this implies  that $L^*(p,x)=\langle x,p\rangle$ and therefore $N(x,p)=\langle x,p\rangle$ and $p\in \bar \partial N (x)$.  

The fact that $N$ is a selfdual Lagrangian follows immediately from the following general lemma. 

\begin{lemma} Let $f_1, f_2: {\cal E}\to \R \cup\{+\infty\}$ be two convex lower semi-continuous functions on a reflexive Banach space ${\cal E}$. The Legendre dual of the function $h$ defined for  $X\in {\cal E}$ by
\[
h(X):=\inf \left\{\frac{1}{2}f_1(X_1)+\frac{1}{2}f_2( X_2)+\frac{1}{8}\|X_1-X_2\|^2; \, X_1, X_2\in {\cal E}, \, X=\frac{1}{2}(X_1+ X_2)\right\}
\]
is given by the function $h^*$ defined for $P\in {\cal E}^*$ by 
\[
h^*(P):=\inf \left\{\frac{1}{2}f^*_1(P_1)+\frac{1}{2}f^*_2( P_2)+\frac{1}{8}\|P_1-P_2\|^2; \, P_1, P_2\in  {\cal E}^*, \, P=\frac{1}{2}(P_1+ P_2)\right\}
\]
\end{lemma} 
{\bf Proof of lemma:} Note that 
\[
h(X):=\inf \left\{F(X_1, X_2); \, X_1, X_2\in {\cal E}, \, X=\frac{1}{2}(X_1+ X_2)\right\}
\]
 where $F$ is the function on $ {\cal E}\times  {\cal E}$ defined as $F(X_1, X_2)= g_1(X_1, X_2)+ g_2(X_1, X_2)$ with 
\[
\hbox{$g_1(X_1, X_2)=\frac{1}{2}f_1(X_1)+\frac{1}{2}f_2( X_2)$\quad  and \quad  $g_2(X_1, X_2)=\frac{1}{8}\|X_1-X_2\|^2$.}
\] 
It follows that  
\[
h^*(P)=F^*(\frac{P}{2}, \frac{P}{2})=(g_1+g_2)^*(\frac{P}{2}, \frac{P}{2})= g_1^*\star g_2^* (\frac{P}{2}, \frac{P}{2}).
\]
It is easy to see that 
\[
g_1^*(P_1, P_2)=\frac{1}{2}f^*_1(\frac{P_1}{2})+\frac{1}{2}f^*_2(\frac{P_2}{2}), 
\]
while 
\[
\hbox{$g_2^*(P_1, P_2)=2\|P_1\|^2$ \quad  if $P_1+P_2=0$\quad  and \quad $+\infty$ otherwise.}
\]
It follows that 
\begin{eqnarray*}
h^*(P)&=&g_1^*\star g_2^* (\frac{P}{2}, \frac{P}{2})\\
&=&\inf \left\{\frac{1}{2}f^*_1(\frac{P_1}{2})+\frac{1}{2}f^*_2( \frac{P_2}{2})+2\|\frac{P}{2}-\frac{P_1}{4}\|^2; \, P_1, P_2\in  {\cal E}^*, \, P=P_1+ P_2\right\}\\
&=&\inf \left\{\frac{1}{2}f^*_1(Q_1)+\frac{1}{2}f^*_2(Q_2)+2\|\frac{P}{2}-\frac{Q_1}{2}\|^2; \, Q_1, Q_2\in  {\cal E}^*, \, P=\frac{1}{2}(Q_1+ Q_2)\right\}\\
&=&\inf \left\{\frac{1}{2}f^*_1(Q_1)+\frac{1}{2}f^*_2(Q_2)+\frac{1}{8}\|Q_2-Q_1\|^2; \, Q_1, Q_2\in  {\cal E}^*, \, P=\frac{1}{2}(Q_1+ Q_2)\right\}.
\end{eqnarray*}

{\bf End of proof of Theorem \ref{ghoussoub}:} Associate to the maximal monotone operator $T$ the ``sub-selfdual" Lagrangian $L_T$ via Lemma \ref{fitz}, that is 
\[
\hbox{$T=\bar \partial L_T$ and $L_T^*(p, x)\geq L_T(x,p)\geq \langle x,p\rangle$.}
\]
Now apply the preceding Proposition to $L_T$ to find a selfdual Lagrangian $N_T$ such that $ L_T(x,p) \leq N_T(x,p) \leq L_T^*(p,x)$ for every $(x,p)\in X\times X^*$, and 
 $Tx= \bar \partial N_T(x)$ for any $x\in  D (T)$.
 
 \section{Operations on maximal monotone operators and the corresponding Lagrangian calculus} 
 
 For a given maximal monotone operator $T: D(T) \subset X \to X^*$, we shall from now on denote by $L_T$ the selfdual Lagrangian on $X \times X^*$ given by Theorem \ref{ghoussoub} in such a way that $\bar \partial L_T=T$. We shall then say that {\it $L_T$ is a selfdual potential for  $T$}.  The following propositions describe the selfdual Lagrangian calculus that parallels the well known calculus developed for monotone operators. Most proofs are good exercises in convex analysis and are left to the interested reader. We also refer to \cite{G2}, and the upcoming \cite{G10}.

 \begin{proposition} Let $T$ be a maximal monotone operator on a reflexive Banach space $X$, and let $L_T$ be its  selfdual potential on $X \times X^*$. Then the following hold:
 
 \begin{enumerate}
  \item If   $\lambda >0$, then the vector field $\lambda \cdot T$ defined by $(\lambda \cdot T)(x)=\lambda T(\frac{x}{\lambda})$ is maximal monotone with selfdual potential given by $(\lambda \cdot L_T)(x,p):=\lambda^2L_T(\frac{x}{\lambda}, \frac{p}{\lambda})$. 
   \item  For $y\in X$ and $q\in X^*$,  the vector field $T^{1, y}$ (resp.,  $T^{2, q}$) given by $T^{1, y}(x)=T(x+y)$ (resp.,  $T^{2, q}(x)=T(x)-q$) is maximal monotone with selfdual potential given by $M_y(x,p)=L_T(x+y,p)-\langle y,p\rangle$ (resp., $N_q(x,p)=L_T(x,p+q)-\langle x, q\rangle$).
        \item If $X$ is a Hilbert space, $U$ is a unitary operator ($UU^*=U^*U=I$) on $X$, then the vector field $T_U$ given by $T_U (x)=U^*T(Ux)$ is maximal monotone with selfdual potential given by   $M(x,p):=L_T(Ux, Up)$.
 
 \item If  $\Lambda:X\to X^{*}$ is any bounded skew-adjoint operator, then the vector field $T+\Lambda $ is a maximal monotone operator with selfdual potential given by $M (x, p)= L_T(x, -\Lambda  x+p)$. 

 \item If  $\Lambda:X\to X^{*}$ is an invertible skew-adjoint operator, then the vector field $\Lambda T^{-1}\Lambda  -\Lambda $  is maximal monotone with selfdual potential given by $M(x, p)= L_T(x+\Lambda ^{-1}p, \Lambda  x)$.
 
  \item   If $\phi$ is a convex  lower semi-continuous function on $X\times Y$ where $X, Y$ are reflexive Banach spaces,  if $A: X\to Y^{*}$ is any  bounded linear operator, and if $J$ is the symplectic operator on $X\times Y$  defined by $J(x, y)=(-y,x)$, then the vector field $\partial \phi + (A^*, A)\circ J$ is maximal monotone on $X\times Y$ with selfdual potential given by 
    $
 L((x,y), (p, q))=\phi (x, y)+\phi^*(A^{*}y+p, -Ax+q)
$.

\end{enumerate}

\end{proposition}

\begin{proposition} {\rm (Direct sums of maximal monotone operators)}

\begin{enumerate}

  \item If $T_i$ is maximal monotone on a reflexive Banach space $X_{i}$ for each $i\in I$, then the vector field $\Pi_{i\in I} T_i$ on $\Pi_{i\in I}X_i$ given by $(\Pi_{i\in I} T_i)((x_i)_i)=\Pi_{i\in I}T_i(x_i)$ is maximal monotone with selfdual potential  $M((x_i)_i, (p_i)_i)=\Sigma_{i\in I} L_{T_i} (x_i, p_i)$. 
      \item If $T_1$ (resp., $T_2$) is a maximal operator on $X$ (resp., $Y$),   then for any bounded linear operator $A: X\to Y^{*}$, the vector field defined on $X\times Y$ by $T= (T_1, T_2)+(A^*, A)\circ J$ is maximal monotone with selfdual potential given by $L ((x,y), (p.q)):=L_{T_1}(x, A^{*}y+p)+L_{T_2}(y, -Ax+q)$.

 \end{enumerate} 
 
  \end{proposition}
  
  \begin{proposition} {\rm (Sums and convolutions)} Let $X$ be a reflexive Banach space $X$. 

   \begin{enumerate}
      \item    If $T$ and $S$ are two maximal monotone operators on $X$ such that  $D(T^{-1})-D(S^{-1})$ contains a neighborhood of the origin in $X^*$, then the vector field $T+S$ is maximal monotone with potential given by 
       \[
 (L_T\oplus L_S)(x,p)=\inf\{L_T(x, r) + L_S(x,p-r); r\in X^{*}\}.
 \]

  \item    If $T$ and $S$ are two maximal monotone operators on $X$ such that $D(T)-D(S)$ contains a neighborhood of the origin in $X$, then the vector field $T\star S$ whose  potential is given by 
   \[
(L_T\star L_S) (x,p)=\inf\{L_T(z, p) + L_S(x-z,p); z\in X\}
\]
is maximal monotone. 
     \end{enumerate}
\end{proposition}
{\bf Proof:}  We only prove 2) as 1) is similar and is left to the reader. It suffices to show that the Lagrangian $L_T\star L_S$ is selfdual. For that fix $(q, y)\in X^{*}\times X$ and write:
 \begin{eqnarray*}
  (L_T\star L_S)^{*} (q,y)
  &=&\sup\{\langle q, x\rangle +  \langle y, p\rangle- L_T(z, p)- L_S(x-z,p); (z, x, p)\in X\times X\times X^{*}\}\\
&=&  \sup\{\langle q, v+z\rangle +  \langle y, p\rangle- L_T(z, p)- L_S(v,p); (z, v, p)\in X\times X\times X^{*}\}\\
  &=&  \sup\{-\phi^*(-z, -v, -p) -\psi^*(z, v, p);\,   (z, v, p)\in X\times X\times X^*\}
 \end{eqnarray*}
where $\phi^*(z, v, p)=\langle q, z\rangle + L_T(-z, -p)$ and $\psi^*(z, v, p)=-\langle y, p\rangle -\langle q, v\rangle +L_S(v,p)$.

Note now that 
 \begin{eqnarray*}
\phi (r, s ,x)&=& \sup\{ \langle r, z \rangle +\langle v, s \rangle+\langle x, p \rangle -\langle q, z \rangle - L_T(-z, -p); (z, v, p)\in X\times X\times X^{*}\}\\
&=& \sup\{ \langle r-q, z \rangle +\langle v, s \rangle+\langle x, p \rangle  - L_T(-z, -p); (z, v, p)\in X\times X\times X^{*}\}\\
&=& \sup\{  \langle v, s \rangle  + L_T^*(q-r, -x); v\in X\}
 \end{eqnarray*}
 which is equal to $+\infty$ whenever $s\neq 0$. 
  Similarly we have 
  \begin{eqnarray*}
\psi (r, s ,x)&=& \sup\{ \langle r, z \rangle +\langle v, s \rangle+\langle x, p \rangle +\langle y, p \rangle +\langle v, q \rangle -L_S(v,p); ; (z, v, p)\in X\times X\times X^{*}\}\\
  &=& \sup\{ \langle r, z \rangle +\langle v, q+s \rangle+\langle x+y, p \rangle    -L_S(v,p); ; (z, v, p)\in X\times X\times X^{*}\}\\
  &=&\sup\{  \langle z, r \rangle  + L_S^*(q+s, x+y); z\in X\}
  \end{eqnarray*} 
   which is equal to $+\infty$ whenever $r\neq 0$. 
  If now $D(T)-D(S)$ contains a neighborhood of the origin in $X$, then we can 
 apply the theorem of Fenchel-Rockafellar \cite{ET} to get 
   \begin{eqnarray*}
  (L_T\star L_S)^{*} (q,y)&=& \sup\{-\phi^*(-z, -v, -p) -\psi^*(z, v, p);\,   (z, v, p)\in X\times X\times X^*\}\\
  &=&\inf \{\phi (r, s ,x)+\psi (r, s ,x); (r,s, x)\in X^*\times X^*\times X\}\\
  &=&\inf \left\{\sup\limits_{v\in X} \{  \langle v, s \rangle  + L_T^*(q-r, -x)\}\right.\\
  &&\left.+\sup\limits_{z\in X}\{  \langle z, r \rangle  + L_S^*(q+s, x+y)\}; (r,s, x)\in X^*\times X^*\times X\right\}\\
   &=&\inf \left\{ L_T^*(q, -x)\} + L_S^*(q, x+y)\}; x\in X\right\}\\
    &=&\inf \left\{ L_T(-x, q)\} + L_S(x+y, q)\}; x\in X\right\}\\
    &=&(L_T\star L_S)(y,q).
 \end{eqnarray*}  
The following way of combining maximal monotone operators has many applications, in particular to partial differential systems and evolution equations. 
\begin{proposition}\label{sum}
Consider $(n+1)$ reflexive Banach spaces $Z$,  $X_1, X_2,...., X_n$,  and bounded linear operators $(A_i, \Gamma_i):Z\to X_i\times X_i^*$ for $i=1, ..., n$, such  that the linear operator $
\Gamma:=(\Gamma_1, \Gamma_2, ..., \Gamma_n): Z\to  \Pi_{i=1}^{n} X_i^*$
is an isomorphism,  
and the following identity holds:
\begin{equation} \label{Stokes.1}
\hbox{$\sum\limits_{i=1}^n\langle A_i z,\Gamma_iz \rangle =0$  for all $z\in Z$.} 
\end{equation}
Let $T_i$ be maximal monotone operators on $X_i$ for $i=1, ..., n$, then the vector field defined on $Z$  by 
\[
z\to Tz:=(T_iA_iz-\Gamma_iz)_{i=1}^n
\]
is maximal monotone with selfdual potential given by the following Lagrangian on $Z\times Z^*$
 \[
L_T(z, p)=\sum\limits_{i=1}^nL_{T_i}(A_iz+p_i, \Gamma_iz). 
\]

\end{proposition}

  \noindent{\bf Proof:} Note that  here $Z$ is put  in duality with the space $X_{1} \oplus X_{2}\oplus ...\oplus X_n$ via the formula
 \begin{equation*}
\langle z, (p_1, p_2,..., p_n) \rangle =\sum\limits_{i=1}^n\langle \Gamma_iz, p_i\rangle, 
\end{equation*}
 where $z\in Z$ and $(p_1, p_2,..., p_n)\in  X_{1} \oplus X_{2}\oplus ...\oplus X_n$. Note also that  (\ref{Stokes.1}) yields that \begin{equation} \label{Stokes.2}
\hbox{$\sum\limits_{i=1}^n\langle A_i y,\Gamma_iz\rangle +\langle A_i z,\Gamma_iy \rangle=0$  for all $y, z\in Z$.} 
\end{equation}
 To show that 
 $L_T$ is  selfdual Lagrangian, fix $((q_1, q_2,..., q_n), y) \in  (X_{1} \oplus X_{2}\oplus ...\oplus X_n) \times Z$  and calculate
\begin{eqnarray*}
L_T^{*} (q,y)=\sup\{\sum\limits_{i=1}^n \langle  \Gamma_iz, q_i\rangle +   \langle \Gamma_iy, p_i\rangle
 \left.-\sum\limits_{i=1}^nL_{T_i}(A_iz+p_i, \Gamma_iz)
; z\in Z, p_{i}\in X_i\right\}. 
\end{eqnarray*}
Setting $x_i= A_iz+p_i\in X_i$, 
we obtain that
\begin{eqnarray*}
L_T^{*} (q,y)&=& \sup\{\sum\limits_{i=1}^n\langle  \Gamma_iz, q_i\rangle   
 +\langle  \Gamma_iy, x_i- A_iz\rangle 
\left.- \sum\limits_{i=1}^nL_{T_i}(x_i, \Gamma_iz):
z\in Z, x_i\in X_i\right\} \\
&=& \sup\{ \sum\limits_{i=1}^n\langle  \Gamma_iz, q_i\rangle 
+\langle  \Gamma_iy, x_i\rangle 
+ \langle A_iy, \Gamma_iz\rangle
\left.-\sum\limits_{i=1}^nL_{T_i}(x_i, \Gamma_iz)
; z\in Z, x_i\in X_i\right\} \\
&=& \sup\{\sum\limits_{i=1}^n\langle  \Gamma_iz, q_i+A_iy\rangle 
+\langle  \Gamma_iy, x_i\rangle  
\left.-\sum\limits_{i=1}^nL_{T_i}(x_i, \Gamma_iz)
; z\in Z, x_i\in X_i\right\}. 
\end{eqnarray*}
Since $Z$ can be identified with $X^*_1\oplus X^*_{2}\oplus ... \oplus X^*_n$ via the correspondence $z\to (\Gamma_1z, \Gamma_2z,..., \Gamma_nz)$, we obtain:
\begin{eqnarray*}
L_T^{*} (q,y)&=&\sup\{\sum\limits_{i=1}^n\langle  z_i, q_i+A_iy\rangle 
+\langle  \Gamma_iy, x_i\rangle  
\left.-\sum\limits_{i=1}^nL_{T_i}(x_i, z_i)
; z_i\in X_i^*, x_i\in X_i\right\}\\ 
&=&\sum\limits_{i=1}^n\sup\left.\{\langle  z_i, q_i+A_iy\rangle 
 +\langle  \Gamma_iy, x_i\rangle  
-L_{T_i}(x_i, z_i)
; z_i\in X_i^*, x_i\in X_i\right\}\\ 
 &=&\sum\limits_{i=1}^nL_{T_i}^{*}( \Gamma_iy, q_i+A_iy ).
\end{eqnarray*} 

 \section{Lifting maximal monotone operators to path spaces}
 
 Let $I$ be any finite time interval that we shall take here --without loss of generality-- to be $[0, 1]$ and let $X$ be a reflexive Banach space.  A   {\it time-dependent --possibly set valued-- monotone map   on $[0, 1]\times X$}  ({\it resp., a time-dependent  convex Lagrangian on $[0, 1]\times X\times X^*$}) is a map $T: [0,1]\times X\to 2^{X^*}$ (resp.,  a function $L: [0,1]\times X\times X^*\to \R \cup \{+\infty\}$) such that :
 \begin{enumerate}
 \item  $T$ (resp., $L$) is measurable with respect to   the  $\sigma$-field  generated by the products of Lebesgue sets in  $[0,1]$ and Borel sets in $X$ (resp., in $X\times X^*$). 
 \item For each $t\in [0,1]$, the map $T_t:=T (t, \cdot )$  is monotone  on $X$ (resp., the Lagrangian $L(t, \cdot, \cdot )$) is convex and lower semi-continuous on $X\times X^*$).
 \end{enumerate}
To each time-dependent Lagrangian $L$ on $[0,1]\times X\times X^*$, one can associate a Lagrangian ${\cal L}$ on the path space $L^{\alpha}_{X}[0,1]\times L^{\beta}_{X^*}[0,1]$ ($\frac{1}{\alpha}+\frac{1}{\beta}=1$)  via the formula:
\[
{\cal L} (u,p):=\int_0^1 L(t, u(t), p (t)) dt.
\]
The Fenchel-Legendre dual of ${\cal L}$ in both variables is then defined for any $(q, v)\in L_{X^*}^{\beta} \times L_{X}^{\alpha}$:
\[
{\cal L}^{*}(q,v) = \sup \left\{ \int_0^1 \big\{ \langle q(t), u(t)\rangle + \langle p(t), v(t)\rangle -
     L(t, u(t),p(t))\big\} dt\ ;  (u,p) \in L_X^{\alpha} \times L_{X^*}^{\beta} \right\}.
\]
It is standard to show that ${\cal L}^{*}(p,u) = \int_0^1 L^{*}(t, p(t),u(t))dt$, which means that if $L$  is a time-dependent selfdual Lagrangian on $[0,1]\times X\times X^*$, then ${\cal L}$ is itself a selfdual Lagrangian on $L^{\alpha}_{X}\times L^{\beta}_{X^*}$. 

The following now follows from the above considerations and Theorem \ref{ghoussoub}. 

\begin{proposition} Let  $T$ be a time-dependent maximal monotone operator on $[0, 1]\times X$, then the function $L_T(t,x,p)=L_{T_t}(x, p)$  is a time-dependent  selfdual Lagrangian $L_T$ on $[0, 1]\times X\times X^*$.

Moreover, 
the  operator $\bar T$ defined on $L^{\alpha}_{X}$ by $\bar T(u(t)_t)=(T(t, u(t)))_t$ is maximal monotone with potential given by 
${\cal L}(u,p) = \int_0^1 L_{T_t}(u(t),p(t))dt$.
\end{proposition}

   We now assume for simplicity that $X$ is a Hilbert space --denoted $H$-- and we consider the space $A_H^2:=\left\{ u:[0,1]\ra H;\dot u\in L_H^2\right\}$
consisting of all absolutely continuous arcs $u:[0,1]\ra H$ equipped with the norm
\begin{eqnarray*}
\| u\|_{A_H^2}={\left\{ \big\|\frac{u(0)+u(1)}{2}\big\|_H^2
+\int_0^1 \|\dot u\|_H^2 \, dt\right\} }^{\frac{1}{2}}.
\end{eqnarray*}
 The space $A_H^2$ can be identified with the product space $H\times L_H^2$, in such a way that its dual $(A_H^2)^*$ can also be identified with $H\times L_H^2$ via the formula
\begin{eqnarray*}
{\braket{u}{(p_1,p_0)}}_{A_H^2,H\times L_H^2}
  =\braket{\frac{u(0)+u(1)}{2}}{p_1} +\int_0^1
   \braket{\dot u(t)}{p_0(t)}\, dt
\end{eqnarray*}
where $u\in A_H^2$ and $(p_1,p_0(t))\in H\times L_H^2$.

\begin{theorem} Let  $T$ be a time-dependent maximal monotone operator on $[0, 1]\times H$, and let $S$ be a maximal monotone operator on $H$, then the operator
\[
{\cal T}u=\left(\dot u +T_tu,  S (u(0)- u(1))+\frac{u(0)+u(1)}{2}\right)
 \]
is maximal monotone on  $A_H^2$, with potential given by the selfdual Lagrangian defined on $A_H^2\times {(A_H^2)}^*=A_H^2\times (H\times L_H^2)$ by
\begin{eqnarray*}
{\cal L}(u,p)=\int_0^1 L_{T_t}\big( u(t)+p_0(t),-\dot u(t)\big)\, dt
  +L_S\big( u(1)-u(0)+p_1,\frac{u(0)+u(1)}{2}\big).
\end{eqnarray*}
Moreover, if we have boundedness conditions of the form 
\begin{equation}
\hbox{ $ \int_0^1 L_{T_t}(x(t),0)\, dt\leq C_1\big(1+\|
x\|_{L^2_H}^2\big)$ for all $x\in L_H^2$,} 
\end{equation}
\begin{equation}
\hbox{  
$L_S(a,0)\leq C_2\big(\| a\|_H^2+1\big)$ for all $ a\in H$, }
 \end{equation}
then  the infimum over $A_H^2$ of the non-negative functional 
\[
I(u)={\cal L}(u,0)=\int_0^1 L_{T_t}\big(u(t),-\dot u(t)\big)\, dt
  +L_S\big( u(1)-u(0),\frac{u(0)+u(1)}{2}\big)
  \]
 is zero and is attained at some $v\in A_H^2$ which  solves the following boundary value problem:
 \begin{eqnarray}
 -\dot{v}(t) &\in& T_t (v(t))\quad {\rm for}\quad   t\in [0,1]\\
 - \frac{v(0)+v(1)}{2}
 & \in& S\big(v(0)- v(1)).
 \end{eqnarray}

\end{theorem} 

\noindent{\bf Proof:} It follows from  Proposition \ref{sum} where we have taken  $Z=A^2_H$ as isomorphic to  the product  $X_1^*\times X_2^*$, where $X_1=X_1^*=L_H^2$ and $X_2=X_2^*=H$, via the map: 
\[
u \in  Z \longmapsto (\Gamma_1u, \Gamma_2u):=\left(-\dot u(t), \frac{u(0)+u(1)}{2}\right)
  \in  L_H^2\times H. 
 \]
The inverse map is
  \[
\big( x,f(t)\big) \in  H\times L_{H}^2\longmapsto x+\frac{1}{2}\left(
  \int_0^t f(s)\, ds-\int_t^1 f(s)\, ds\right) \in Z.
\]
 Define now the maps $A_1:Z\to X_1$ by $A_1u=u$ and $A_2Z\to X_2$ by $A_2(u)=u(1)-u(0)$ in such a way that 
\begin{equation}\label{stokes.2}
\langle A_1u, \Gamma_1u\rangle + \langle A_2u, \Gamma_2u\rangle =-\int_0^1 u(t)\dot u(t) dt+\langle u(1)-u(0), \frac{u(0)+u(1)}{2}\rangle=0.
\end{equation}
If the above  boundedness conditions are satisfied, 
then we can use Proposition \ref{main.existence} to conclude that the infimum over $A_H^2$ of the non-negative functional 
$I(u)={\cal L}(u,0)$
 is zero and is attained at some $v\in A_H^2$. 
By writing  that $I(v)=0$, we get 
\begin{eqnarray*}
0 &=&\int_0^1 \left[L_{T_t}\big( v(t),-\dot{v}(t)\big) + \braket
{v(t)}{\dot{v}(t)}\right]\, dt- \int_0^1 \braket { v(t)}{\dot{v}(t)} \,dt + L_S\big(v(1)-v(0),\frac{v(0)+v(1)}{2}\big) \\
&=&  \int_0^1 \left[L_{T_t}\big( v(t),-\dot{v}(t)\big) + \braket
{v(t)}{\dot{v}(t)}\right]\, dt -\frac{1}{2}\langle v(1)-v(0), v(0)+v(1)\rangle+ L_S
  \big(v(1)-v(0),\frac{v(0)+v(1)}{2}\big).
\end{eqnarray*}
Since $L_{T_t}$ and $L_S$ are selfdual Lagrangians, we have
$
 L_{Tt}\big( v(t),-\dot{v}(t)\big) + \braket
{v(t)} {\dot{v}(t)}\geq 0
$
and
\begin{eqnarray*}
 L_S
  \big(v(1)-v(0),\frac{v(0)+v(1)}{2}\big)-  \braket {v(1)- v(0)}{\frac{v(0)+v(1)}{2}} \geq 0.
\end{eqnarray*}
 which means that 
$
 L_{T_t}\big( v(t),-\dot{v}(t)\big) + \braket
{v(t)} {\dot{v}(t)}\,\big) = 0
$ for almost all $t\in [0,1]$, 
and
\begin{eqnarray*}
 L_S\big(v(1)-v(0),\frac{v(0)+v(1)}{2}\big)-  \braket {v(1)- v(0)}{\frac{v(0)+v(1)}{2}} = 0.
\end{eqnarray*}
The result now follows from the above identities and the limiting case in Fenchel-Legendre duality. 
 In other words $v$ solves the following boundary value problem:
 \begin{eqnarray}
 -\dot{v}(t)&\in&T_t\big(v(t)\big)\quad {\rm for}\quad   t\in [0,1]\\
  -\frac{v(0)+v(1)}{2}
 & \in& S\big(v(0)- v(1)\big).
\end{eqnarray}

 \subsection{Semi-groups of contractions associated to maximal monotone operators}
 
We now use the above representation to recover variationally a well known fact:  To any maximal monotone operator $T$ on a Hilbert space $H$, one can associate a semi-group of contractions $(S_t)_{t\in \R^+}$ on $D(T)$ in such a way that for any  $x_0\in D(T)$ $S_tx_0=x(t)$ where $x(t)$ is the unique solution of the equation:
  \begin{eqnarray}
 \left\{ \begin{array}{l}
\label{eqn:0}
\hbox{$ -{\dot x}(t) \in T(x(t))$ for $t\in [0, 1]$}\\
\quad x(0)=x_0. 
\end{array}\right.
\end{eqnarray}
As noted above, such a solution can be obtained by minimizing the non-negative functional
\[
   I(u)=  \int_0^1   L_T(u(t), -{\dot u}(t))dt +\frac{1}{2}|u(0)|^{2} -2\langle x_0, u(0)\rangle +|x_0|^{2} +\frac{1}{2}|u(1)|^{2}
\]
on $A^2_H$ and by showing that $  I(x)=\inf\limits_{u\in A^2_H} I(u)=0$. Note that we have used for a maximal monotone operator on the boundary the subdifferential $S=\partial \phi$ where $\phi$ is the convex function on $H$ given by 
$\phi (x)=  \frac{1}{4}|x|^2-\langle x, x_0\rangle$.

Now according to the preceeding section, this can be done whenever the Lagrangian $L$ satisfies the following  boundedness condition:
\begin{equation}
L(x,0)\le C(\| x\|^2+1) \quad \hbox{\rm for all $x\in H$, }
\end{equation}
which is too stringent. This condition can however be weakened  by using a $\la$-regularization procedure of selfdual Lagrangians.  The details are given in \cite{GT3} where the following more general result is established. 
  \begin{theorem} \label{semi-group}Let  $T$ be a maximal monotone operator on a Hilbert space $H$ such  $D(T)$ is non-empty, and let $L_T$ be a corresponding selfdual potential.  Then for any $\omega \in {\bf R}$, there exists a semi-group of  maps $(S_{t})_{t\in {\bf R}^{+}}$ defined on $D(T)$  such that: 
\begin{enumerate} 
\item   $S_{0}x=x$ and $\|S_tx-S_ty\|\leq e^{-\omega t}\|x-y\|$ for any $x,y\in D(T)$. 

\item For any  $x_0\in D(T)$, we have $S_tx_0=x(t)$ where $x(t)$ is the unique path that minimizes  on $A_{H}^2$ the following  functional   
\[
   I(u)=  \int_0^1 e^{2\omega t} L_T(u(t), -\omega u (t)- {\dot u}(t))dt +\frac{1}{2}\|u(0)\|^{2} -2\langle x_0, u(0)\rangle +\|x_0\|^{2} +\frac{1}{2}\|e^{\omega}u(1)\|^{2}.
\]

\item For any $x_0\in D(T)$ the path $x(t)=S_tx_0$ satisfies:
 \begin{eqnarray}
\label{eqn:2.2}
 -{\dot x}(t)-\omega x(t)&\in& T(x(t)) \quad \hbox{ \rm for $t\in [0, 1]$}\\
x(0)&=&x_0. \nonumber
\end{eqnarray}
 \end{enumerate}
 \end{theorem}

\subsection{Connecting two maximal monotone graphs}

Let now  $Z=A^2_H\times A^2_H$ with $H$ being a Hilbert space. We shall identify  it with the product space $X^*_1\oplus X^*_2 \oplus X^*_3$, with  $X_1=X^*_1=L^2_{H}\times L^2_{H}$, $X_2=X_2^*= H$ and $X_3=X_3^*= H$ in the following way:
\[
(u,v) \in Z \longmapsto (\Gamma_1(u,v), \Gamma_2(u,v)), \Gamma_3(u,v)=\big((\dot u(t), \dot v(t)), u(0), v(T) \big)  \in (L^2_{H}\times L_H^2)\times H\times H.
\]
The inverse map is then
\[
\big( (f(t), g(t)), x, y\big) \in  (L^2_{H}\times L_H^2)\times H\times H \longmapsto (x+ \int_0^tf(s)ds, \, y-\int_t^1g(s)ds) \in Z=A^2_H\times A^2_H.
\]
The map $A_1:Z\to X_1:=L^2_{H}\times L^2_{H}$ is defined as $A_1(u,v)=(v,u)$ while $A_2:Z\to X_2:=H$ is defined as $A_2(u,v)=v(0)$ and $A_3(u,v)=-u(1)$. It is clear that 
\begin{equation}
\sum\limits_{i=1}^3\langle A_i(u,v), \Gamma_i(u,v)\rangle 
=\int_0^1(\dot u (t) v(t) +\dot v (t) u(t) ) dt+\langle u(0), v(0)\rangle -\langle u(1), v(1)\rangle=0.
\end{equation}
The following now follows readily from Theorem \ref{ghoussoub}, and Propositions \ref{sum} and \ref{main.existence}. 

\begin{theorem} \label{connecting.2}
Suppose $(T_t)_t$ is a time-dependent maximal monotone operator on $[0,1]\times H^2\times H^2$ and that $S_1$, $S_2$ are  two maximal monotone operators  on $H$. The following map 
\[
(u, v) \to\left( - (\dot v, \dot u) +T_t (u, v),  -u(0)+S_1v(0), u(1)+S_2v(1) \right)
\]  
is then maximal monotone on the space $A_H^2\times A_H^2$ whose dual had been identified with $L_{H^2}^2\times H\times H$. The corresponding selfdual potential is given  for $(U, P):=\big((u,v), (p^1_0(t), p^2_0(t), p_1, p_2)\big)\in  Z\times Z^*=(A_H^2\times A_H^2)\times  (L_{H^2}^2\times H\times H)$ by
\begin{eqnarray*}
{\cal L}(U, P)=\int_0^1 L_{T_t}\big((v+p^1_0,u+p^2_0), (\dot u, \dot v)\big)\, dt
  +L_{S_1}( v(0)+p_1, u(0)) +L_{S_2}( -u(1)+p_2, v(1)).
\end{eqnarray*}
  If in addition   we have the following boundedness conditions:
\begin{equation}
\hbox{ $ \int_0^1 L_{T_t}(q_1(t),q_2(t)), 0)\, dt\leq C\big(1+\|
q_1\|_{L^2_H}^2+\|q_2\|_{L^2_H}^2\big)$ for all $(q_1, q_2)\in L_{H^2}^2$,} 
\end{equation}
\begin{equation}
\hbox{ $L_{S_i} (a,0)\leq C_i\big(\| a\|_H^2+1\big)$ for all $ a\in H$, }
 \end{equation}
then the infimum over $Z:=A_H^2\times A_H^2$ of the non-negative functional    
\[
I(u,v)={\cal L}((u(t),v(t)),0)=\int_0^1 L_{T_t}\big( (v(t), u(t)), (\dot u(t), \dot v(t)\big)\, dt
  +L_{S_1}\big(v(0), u(0)) \big)+L_{S_2}\big(-u(1), v(1)) \big)
  \]
is zero and is attained at some $(\bar u, \bar v)\in A_H^2\times A_H^2$ which solves the following boundary value problem:
 \begin{eqnarray}
 (\dot{v}(t) , \dot{u}(t)) &\in& T_t\big(u(t), v(t))\quad {\rm for}\quad   t\in [0,1]\\
  u(0) &\in&  S_1(v(0))\\
  - u(1) &\in& S_2(v(1)).
   \end{eqnarray}
 
\end{theorem}
  \section{Variational resolution and inverse problems of equations involving maximal monotone operators}
 
 In this section, we shall illustrate  how one can use the Lagrangians associated with a maximal monotone operators to give a variational resolution for certain non-potential equations as well as a minimizing procedures for corresponding inverse problems. We limit ourselves here to 
 the following standard elliptic problem
\begin{equation}
\label{Ex50}
 \left\{ \begin{array}{lcl}
    \hfill -{\rm div} (T\nabla u(x))+\lambda u (x)&=& g(x) \quad  \hbox{\rm \, on \, $\Omega \subset \R^n$}
\\
 \hfill  u(x) &=& 0 \quad \quad   \hbox{\rm on \quad $\partial \Omega$, }  
       \end{array}  \right.
   \end{equation}
and to the corresponding dynamic problem below (\ref{Ex51}). Here $g\in L^2(\Omega)$, $\lambda \in \R$, and $T: \R^n \to \R^n$ is a given maximal monotone mapping.

\subsection{A general variational principle for selfdual Lagrangians}

The following general variational principle  can be seen as an extension of Proposition \ref{main.existence}. We recall that the co-Hamiltonian $\tilde H_L$ associated to a Lagrangian $L$ is defined as the Legendre transform of $L$ with respect to the first variable, that is for $(p, q)\in X^*\times X^*$ we have
\[
\tilde H_L(p,q)=\sup\{\langle y,p\rangle -L(y,q); y\in X\}.
\]
 It is  easy to see that  if $L$ is a  selfdual Lagrangian on $X\times X^*$, then its  co-Hamiltonian ${\tilde H}_L$ on $X^*\times X^*$ satisfies the following properties: 
   \begin{itemize}
\item  For  each $p\in X^*$, the function $q\to {\tilde H}_L(p,q)$\big)
 is concave. 

\item  For each $q\in X^*$, the function $p\to {\tilde H}_L(p,q)$\big) is convex and lower semi-continuous.

\item For  each $p, q\in X^*$, we have  ${\tilde H}_L(q,p)\leq -{\tilde H}_L(p,q)$\big) and therefore  ${\tilde H}_L(p,p)\leq 0$. 
\end{itemize}
For more details, we refer to  \cite{G3} and \cite{G10}. In the sequel, we shall use the following min-max theorem of Ky-Fan \cite{AE}. 

\begin{lemma} \label{KF}Let  $E$ be  a  closed convex  subset of a reflexive Banach space $Z$, and consider  $M :E\times E\ra \bar \R $ to be a functional such that
\begin{enumerate}
\item   For each $y\in E$, the map $x \mapsto M(x,y)$ is weakly lower semi-continuous on $E$;
\item   For each $x\in E$, the map $y\mapsto M(x,y)$ is  concave on E;
\item There exists $\gamma \in \R$ such that $M(x,x) \leq \gamma$ for every $x\in E$; 
\item There exists a   $y_0\in E$ such that $E_0= \{ x\in E: M (x,y_0)\leq\gamma\}$ is bounded. 
\end{enumerate}
Then there exists $\bar x \in E$ such that $M(\bar x,y)\leq\gamma$ for all $y\in E$.

\end{lemma} 
We can now prove the following  result. 
 
\begin{proposition} \label{superposition.2} Consider $n$ reflexive Banach spaces $X_1, X_2,...., X_n$,  and let $L_i$ be selfdual  Lagrangians on $X_i\times X_i^*$ for $i=1,...n$.  Let $Z$ be a reflexive Banach space and let $\Gamma_i:Z\to X_i^*$ be bounded linear operators  such that $\Gamma:=(\Gamma_1, \Gamma_2, ..., \Gamma_n)$ maps a closed convex subset $E$ of $Z$ onto $\Pi_{i=1}^nX_i^*$.
Consider also bounded linear operators  $A_i:Z\to X_i$ such that  the map
\begin{equation}
 \hbox{\rm $z\to \sum\limits_{i=1}^n\langle A_i z,\Gamma_i  z \rangle$ is weakly upper semi-continuous on $E$.}
 \end{equation}
Assume that the following coercivity condition holds: 
 \begin{equation}
\hbox{$\lim\limits_{z\in E, \|z\|\to +\infty}\sum\limits_{i=1}^n{\tilde H}_{L_i}(\Gamma_iz, 0) -\langle A_i z,\Gamma_i z \rangle=+\infty$.}
\end{equation} 
Then the infimum of the functional
\begin{equation}
I(z)=\sum\limits_{i=1}^nL_i(A_iz, \Gamma_iz)-\langle A_i z,\Gamma_i z \rangle
\end{equation}
over $E$ is zero and is attained at some $\bar z \in E$,   which then  solves the system of equations:
 \begin{equation}
  \label{NG.1}
 \Gamma_i \bar z_i \in \bar\partial L_i(A_i\bar z_i)\quad {\rm for} \quad i=1,..., n.
  \end{equation}
 \end{proposition}
\noindent{\bf Proof:} First, it is clear that $I$ is non-negative since each $L_i$ is selfdual and therefore  $L_i(x, p) \geq \langle x, p\rangle$ for every $i=1,...,n$. 
Now write for any $z\in Z$
\begin{eqnarray*}
I(z)&=&\sum\limits_{i=1}^nL_i(A_iz, \Gamma_iz)-\langle A_i z,\Gamma_i z \rangle\\
&=&\sum\limits_{i=1}^nL_i^*(\Gamma_i z, A_i z) -\langle A_i z,\Gamma_i z \rangle\\
 &=&\sum\limits_{i=1}^n\sup\big\{ \langle \Gamma_iz,x_i\rangle) +\langle  A_i z, w_i\rangle -L_i(x_i, w_i); x_i\in X_i, w_i\in X_i^*\big\} -\sum\limits_{i=1}^n\langle A_i z,\Gamma_i z \rangle \\
&=&\sum\limits_{i=1}^n\sup\big\{ {\tilde H}_{L_i}(\Gamma_i z, w_i)+\langle A_i z, w_i\rangle;  w_i\in X_i^*\big\} -\sum\limits_{i=1}^n\langle A_i z,\Gamma_i z \rangle \\
&=&\sup\left\{ \sum\limits_{i=1}^n{\tilde H}_{L_i}(\Gamma_i z, \Gamma_iw)+\langle A_i z, \Gamma_i w\rangle;  w\in E\right\} -\sum\limits_{i=1}^n\langle A_i z,\Gamma_i z \rangle \\
&=&\sup\left\{\sum\limits_{i=1}^n {\tilde H}_{L_i}(\Gamma_i z, \Gamma_iw)+\langle A_i z, \Gamma_i (w-z)\rangle;  w\in E\right\}.
 \end{eqnarray*}
In other words,
\begin{equation}
I(z)=\sup\limits_{w\in E}M(z, w)
\end{equation}
where
 \begin{equation}
 M(z,w)=\sum\limits_{i=1}^n{\tilde H}_{L_i}(\Gamma_iz, \Gamma_iw) +\langle A_i z,\Gamma_i (w-z) \rangle,
 \end{equation}
and  where for each $i=1,..., n$, ${\tilde H}_{L_i}$ is the co-Hamiltonian on $X_i^*\times X_i^*$ associated to $L_i$. 

Note now that for each $w\in E$, the map $z \mapsto M(z,w)$ is weakly lower semi-continuous on $E$, and 
for each $z\in E$, the map $w\mapsto M(z,w)$ is  concave on E. Moreover,  $M(z,z) \leq 0$ for every $z\in E$. 
Finally,  use the coercivity condition and Lemma \ref{KF} to deduce the existence of $\bar z$ such that
\[
0\geq \sup\limits_{w\in E}M(\bar z, w)=I(\bar z)= \sum\limits_{i=1}^nL_i(A_i\bar z, \Gamma_i\bar z)-\langle B_iA_i \bar z,\Gamma_i \bar z \rangle\geq 0.
\]
Again, since each term is non-negative, we get that $L_i(A_i\bar z, \Gamma_i\bar z)-\langle A_i \bar z,\Gamma_i \bar z \rangle=0$ for each $i=1,...,n$, and we are done.

\subsection{Variational resolution for non-potential equations and evolutions}
In order to resolve (\ref{Ex50}), we use Theorem \ref{ghoussoub} to associate to the maximal monotone operator $T^{-1}$, a selfdual Lagrangian $L$ on $\R^n\times \R^n$ such that $\bar \partial L=T^{-1}$. We claim that the infimum of the  functional
 \begin{equation}
J(u)=\int_\Omega \left\{L_T\big(\nabla (-\Delta)^{-1}(-\lambda u +g)(x), \nabla u(x)\big) +\lambda |u(x)|^2-   u(x) g(x)  \right\}dx
 \end{equation}
on $H^1_0(\Omega)$ is zero and is attained  at some $u$ which is a solution of (\ref{Ex50}). 

Indeed, we first define the selfdual Lagrangian ${\cal L}_T$ on $L^2(\Omega; \R^n)\times L^2(\Omega; \R^n)$ via the formula
\[
{\cal L}_T(u,p)=\int_\Omega L(u(x), p(x)) dx. 
\] 
 Consider now the space $Z=H^1_0(\Omega)$ and let $X^*$ be the closed subspace of $L^2(\Omega; \R^n)$ defined as 
 \[
\hbox{$ X^*=\{f\in L^2(\Omega; \R^n); \nabla u=f$ for some $u\in H^1_0(\Omega)\}.$}
 \]
 It suffices now to apply Proposition \ref{superposition.2} with the selfdual Lagrangian ${\cal L}_T$, and the operators $A: Z\to X$, and $\Gamma: Z\to X^*$  defined by 
 \[
\hbox{$Au=\nabla (-\Delta)^{-1}(-\lambda u +g)$ \quad and  \quad $\Gamma u =\nabla u$ respectively. }
 \]
 Note that $\Gamma$ is onto and that the diagonal map 
\[
u\to \langle Au, \Gamma u \rangle =\int_\Omega \left\langle \nabla (-\Delta)^{-1}(-\lambda u +g)(x), \nabla u(x)\right\rangle dx = -\lambda \int_\Omega |u|^2\, dx+\int_\Omega g u\, dx 
\]
 is clearly weakly  continuous on $Z$. One also needs a coercivity condition of the type
 \begin{equation}
 L_T(y, 0) \leq C(1+ |y|^2)
 \end{equation}
  since then 
\begin{eqnarray*}
 {\tilde H}_{{\cal L}_T}(\Gamma u, 0) -\langle A u,\Gamma u \rangle&=&  \int_\Omega {\tilde H}_{L_T}(\nabla u, 0) dx +\lambda \int_\Omega |u|^2\, dx-\int_\Omega g u\, dx\\
 &\geq& \frac{1}{C}\int_\Omega |\nabla u|^2 dx +\lambda \int_\Omega |u|^2\, dx-\int_\Omega g u\, dx-C\\
&\geq& 
K_1\int_\Omega |\nabla u|^2\, dx-K_2(\int_\Omega |\nabla u|^2\, dx)^{1/2} -C 
 \end{eqnarray*}
which goes to $+\infty$ with $\|u\|_{H^1_0}$ provided $\lambda >  -\frac{\lambda_1}{C}$.
 
 Let now $v\in H^1_0(\Omega)$ be such that  
 \[
 J(v)=\int_\Omega \left\{L_T\big(\nabla (-\Delta)^{-1}(-\lambda v +g), \nabla v\big) -\lambda |v|^2+\langle v, g \rangle \right\}dx
=0.
\]
It follows that $\nabla (-\Delta)^{-1}(-\lambda v +g) \in (\bar \partial {\cal L}_T)^{-1} (\nabla v)$ and that $-\lambda v +g =-{\rm div} \big(\bar \partial {\cal L}_T)^{-1} (\nabla v)\big)$. In other words, $-\lambda v +g =-{\rm div} (T(\nabla v))$.\\

Similarly, one can solve the corresponding dynamic problem, 
\begin{equation}
\label{Ex51}
 \left\{ \begin{array}{lcl}
    \hfill u_t(t, x)-{\rm div} (T\nabla_x u(t, x))&=& g(t, x) \quad \, \,  \hbox{\rm \, on \, $[0, 1]\times \Omega$}
\\
 \hfill  u(x) &=& 0 \quad \quad  \quad  \quad \hbox{\rm on \quad $\partial \Omega$, }  \\
\hfill  \frac{u(0)+u(1)}{2}
 & =& S\big(u(1 )- u(0)), 
       \end{array}  \right.
   \end{equation}
where $g\in L^2(\Omega)$, $T: \R^n \to {\R^n}$ is a given maximal monotone mapping on $\R^n$ and $S$ is a maximal monotone operator on $L^2(\Omega)$. It suffices to minimize the non-negative functional 
 \begin{equation}
{\cal J}(u)=\int_0^1\int_\Omega \left\{L_T\big( \nabla (-\Delta)^{-1}(-\dot u +g), \nabla u\big) -  u g  \right\}dx dt + L_S\big( u(1)-u(0),\frac{u(0)+u(1)}{2}\big)
\end{equation}
on  $W^{1,2}([0, 1]; H^1_0(\Omega))$ where $L_S$ is the selfdual potential  on $L^2(\Omega)\times L^2(\Omega)$ associated to $S$.

Indeed, on can apply Proposition \ref{superposition.2} with
\begin{itemize}

\item ${\cal Z}:=W^{1,2}([0, 1]; H^1_0(\Omega))$, ${\cal X}_1:=L^2([0,1]; X)$ and ${\cal X}_2:=L^2(\Omega)$, where  the space $X^*$ is again given by  $X^* =\{f\in L^2(\Omega; \R^n); \nabla u=f$ for some $u\in H^1_0(\Omega)\}.$
 
\item the operators $A_1: {\cal Z}\to {\cal X}_1=L^2([0, 1];X)$ resp., $A_2: {\cal Z}\to {\cal X}_2:=L^2(\Omega)$ are given by 
\[
\hbox{$A_1u=\nabla (-\Delta)^{-1}(-\dot u +g)$ \quad and  \quad $A_2 u(x) = u(1, x)- u(0, x)$ respectively, }
 \]
while $\Gamma_1: {\cal Z}\to {\cal X}^*_1=L^2([0, 1];X^*)$, and $\Gamma_2: {\cal Z}\to {\cal X}^*_2=L^2(\Omega)$ are defined by 
 \[
\hbox{$\Gamma_1 u =\nabla_x u$ \quad and \quad  $\Gamma_2u(x)=\frac{u(0, x )+u(1, x )}{2}$ respectively. }
 \]
Note that $(\Gamma_1, \Gamma_2):Z\to {\cal X}^*_1\times {\cal X}^*_2$ is an isomorphism and that the diagonal map 
\begin{eqnarray*}
\langle A_1u, \Gamma_1u\rangle + \langle A_2u, \Gamma_2u\rangle &=&-\int_0^1 \langle u(t), \dot u(t)-g(t)\rangle_{L^2(\Omega)} dt+\langle u(1)-u(0), \frac{u(0)+u(1)}{2}\rangle\rangle_{L^2(\Omega)}\\
&=&\int_0^1\int_\Omega u(t,x)g(t, x)\, dxdt
\end{eqnarray*}
is clearly weakly continuous on $Z$.
 
\item  the selfdual Lagrangians ${\cal L}_T$, $L_S$ on ${\cal X}_1:=L^2([0,1]; X)$ and ${\cal X}_2:=L^2(\Omega)$ respectively,  where ${\cal L}_T$ is given by 
${\cal L}_T (U, P)=\int_0^1\int_\Omega L_T(U(t, x), P(t, x)) dx dt$.
\end{itemize}
As above, if we have  boundedness conditions of the form 
\begin{equation*}
\hbox{ $ L_{T}(x,0)\, dt\leq C_1\big(1+|
x|^2\big)$ for $x\in \R^n$ \quad and \quad
 $L_S(u,0)\leq C_2\big(\| u\|_2^2+1\big)$ for $ u\in L^2(\Omega)$, }
 \end{equation*}
then  the infimum over $Z$ of the non-negative functional ${\cal J}$ is zero and is attained at some $v\in Z$ which  solves the  
evolution equation (\ref{Ex51})

 \subsection{Inverse problems}
 
 We consider here the following inverse problem: {\it Given $u_0\in H^1_0(\Omega)$, find a maximal monotone vector field $T$ in a given class ${\cal C}$ such that $u_0$ is a solution of the corresponding equation (\ref{Ex50}).}

Since such a vector field $T$ may not exist in general in the class ${\cal C}$, so  one proceeds to find a maximal monotone vector field $T \in {\cal C}$, such that the corresponding solution $u$ of  (\ref{Ex50}) is as close as possible to $u_0$. The least square approach leads to the following minimization problem
\begin{equation}\label{min}
\hbox{$\inf \{\int_\Omega |u(x)-u_0(x)|^2 dx;\,  u\in H^1_0(\Omega), T \in {\cal C},$ such that $  -{\rm div} (T(\nabla u))+\lambda u=g$  on $\Omega$ \}.}
\end{equation}
 Note that the constraint set above is not easily tractable, but in view of our variational characterization of the solutions $u$, one is able to approach the problem via the following penalized least square minimization procedure. 
 
 Let ${\cal L}$ be the class of selfdual Lagrangians corresponding to ${\cal C}$, that is  
\[
\hbox{${\cal L}=\{ L$ selfdual on $\R^n\times \R^n;\,  \bar \partial L=T^{-1}$ for some $T\in {\cal C}$\}.}
\]
 For each $\epsilon>0$, we consider  the minimization problem:
  \begin{equation}\label{min.epsilon}
\inf \{ {\cal P}_\epsilon (L, u);\,  L\in {\cal L}, u\in H^1_0(\Omega)\}.
  \end{equation}
where 
\[
{\cal P}_\epsilon (L, u)= \int_\Omega |u(x)-u_0(x)|^2 dx +\frac{1}{\epsilon}\int_\Omega \left\{L\big(\nabla (-\Delta)^{-1}(-\lambda u +g), \nabla u \big) +\lambda |u|^2-   u g  \right\}dx.
\]
Note that ${\cal P}_\epsilon$ is convex and lower semi-continuous in both variables $(L, u)$, and therefore if ${\cal L}$ is a suitable convex compact class of selfdual Lagrangians, then there exists  a minimizer $(L_\epsilon, u_\epsilon)\in {\cal L}\times H^1_0(\Omega; \R^n)$ for (\ref{min.epsilon}). 
Now when $\epsilon$ is small enough, 
 the non-negative penalization has to be very small at the minimum $(L_\epsilon, u_\epsilon)$.  
 In other words, any weak cluster point $(L_0, u_0)$ of the family $(L_\epsilon, u_\epsilon)_\epsilon$  is a solution of problem (\ref{min}) with $T_0:=\bar \partial L_0$ being the optimal maximal monotone operator, since the penalty term 
$\int_\Omega \left\{L\big(\nabla (-\Delta)^{-1}(-\lambda u +g), \nabla u\big) +\lambda |u|^2-   u g  \right\}dx$  has to be zero.   

For more details, we refer the reader  to \cite{Z}.


\begin{thebibliography}{99}

 \bibitem{AE} J.P. Aubin, I. Ekeland {\em Applied nonlinear analysis},  Reprint of the 1984 original. Dover Publications, Inc., Mineola, NY, 2006.

\bibitem{BW} H.H. Bauschke, X. Wang, {\em The kernel average of two convex functions and its applications to the extension and representation of monotone operators}, Preprint (2007).

 \bibitem{Br} H. Brezis,  {\em Op\'erateurs maximaux monotones et semi-groupes
de contractions dans les espaces de Hilbert},  North Holland, 
Amsterdam-London (1973).

 \bibitem{Bro1} F. Browder, {\em Probl\`emes non lin\'eaires}, Presses de l'Universit\'e de Montr\'eal (1966)

 \bibitem{Bro2} F. Browder, {\em Nonlinear maximal monotone operators in Banach space}, Math. Annalen {\bf 175} (1968) p. 89-113.
 
\bibitem{BS} R.S. Burachek, B. F. Svaiter, {\em Maximal monotonicity, conjugation and the duality product}. P.A.M.S, {\bf 131}, 8 (2003) p. 2379-2383.  


\bibitem{BE} H. Brezis, I. Ekeland, {\em Un principe variationnel
  associ\'e \`a certaines equations paraboliques. Le cas independant du
 temps}, C.R. Acad. Sci. Paris S\'er. A {\bf 282} (1976), 971--974.


\bibitem{F} S.P. Fitzpatrick, {\em Representing monotone operators by convex functions}, Proc. Centre for Math. Analysis 20 (1989), p. 59-65. 

\bibitem{ET} I. Ekeland, R. Temam, {\em Convex Analysis and Variational problems},
Classics in Applied Mathematics, {\bf 28} SIAM (1999 Edition).

\bibitem{Ev} L.C. Evans,  {\em Partial Differential Equations},Graduate Studies in Mathematics, vol. 19,
Amer. Math. Soc., Providence, 1998.

 \bibitem{G2} N. Ghoussoub,  {\it Anti-selfdual Lagrangians: Variational resolutions of  non self-adjoint equations and dissipative evolutions},  AIHP-Analyse non lin\'eaire, 24 (2007) p. 171-205.
 
 \bibitem{G3} N. Ghoussoub, {\it Anti-symmetric Hamiltonians: Variational resolution of Navier-Stokes equations and other nonlinear evolutions},  Comm. Pure \& Applied Math., vol. 60, no. 5 (2007) pp. 619-653.
 
 \bibitem{G4}  N. Ghoussoub,  {\it Superposition of selfdual functionals for non-homogeneous boundary value problems and differential systems}, Journal of Discrete and Continuous Dynamical Systems, Vol. 21, 1 (2008), p. 71-104.  

 \bibitem{G5} N. Ghoussoub, {\it  Maximal monotone operators are selfdual vector fields and vice-versa},   arXiv:math/0610494  (2006) 9 pages.
 
 \bibitem{G6}   N. Ghoussoub,   {\it  Hamiltonian systems as selfdual equations}, Frontiers in Mathematics in China,   (2008) 19 pp.

\bibitem{G10} N. Ghoussoub, {\it  Selfdual partial differential systems and their variational principles},  Springer-Verlag, Universitext (2007) (To appear) 350 pp.

\bibitem{GM2} N. Ghoussoub, A. Moameni, {\it   Selfdual variational principles for periodic solutions of Hamiltonian and other dynamical systems},   Comm. in PDE  32, (2007) p. 771-795.

\bibitem{GM3} N. Ghoussoub, A. Moameni, {\it   Hamiltonian systems of PDEs with selfdual  boundary conditions}, submitted  (2007)  30 p.

\bibitem{GM4}  N. Ghoussoub, A. Moameni: {\it Anti-symmetric Hamiltonians (II): Variational resolution of Navier-Stokes equations and other nonlinear evolutions}, Annales I.H.P, Analyse non-lin\'eaire, In press (Accepted July 2007) 26 p.

\bibitem{GT1} N. Ghoussoub, L. Tzou,  {\em A variational principle for gradient flows}, Math. Annalen, Vol 30, 3 (2004) p. 519-549. 

\bibitem{GT2} N. Ghoussoub, L. Tzou,  {\it Anti-selfdual Lagrangians II: Unbounded non self-adjoint operators and  evolution equations},  Annali di Matematica Pura ed ApplicataÓ, Vol 187, 2 (2008) p. 323-352 (Published online March 30, 2007) 


\bibitem{GT3} N. Ghoussoub, L. Tzou,   {\em  Iterations of anti-selfdual Lagrangians and applications to Hamiltonian systems and multiparameter gradient flows}, Calc. Var.\& PDE  Vol. 26, N. 4 (2006) p. 511- 534 

  \bibitem{KS} D. Kinderlehrer, G. Stampachia, {\em An introduction to variational inequalities and their applications}, Classics in Applied Math, 31, SIAM (2000) 

\bibitem{K} E. Krauss, {\em A representation of maximal monotone operators by saddle functions}, Rev. Roum. Math. Pures Appl. 309 (1985) p. 823-837

\bibitem{Ph} R.R Phelps, {\em Convex functions, monotone operators and
differentiability}, Lecture Notes in Math. 1364, Springer Verlag, New York, Berlin, Tokyo,
(1998), 2nd edition 1993.

\bibitem{S} B. F. Svaiter, {\em Fixed points in the family of convex representations of a maximal monotone operator,}  P.A.M.S, {\bf 131}, 12 (2003) p. 3851-3859.  

% \bibitem{Te2} R. Temam, {\it Infinite-dimensional dynamical systems in mechanics and physics}, Applied mathematical sciences, 68, Springer-Verlag (1997).


\bibitem{Z} R. Zarate, PhD dissertation, The University of British Columbia (2008). 


\end{thebibliography}
\end{document}